\documentclass[oneside,11pt]{article}

\usepackage{amsmath}
\usepackage{amssymb}
\usepackage{pxfonts}
\usepackage{graphicx}
\usepackage{theorem}
\usepackage{pifont}
\usepackage{color}
\usepackage[all]{xy}
\usepackage{tikz}
\usepackage{tkz-euclide}

\usetkzobj{all}

\newcommand{\beq}{\begin{equation}}
\newcommand{\eeq}{\end{equation}}
\newcommand{\bbR}{{\mathbb R}}
\providecommand{\R}{\bbR}

\newcommand{\calM}{{\mathcal M}}

\newcommand{\Textand}{\qquad\text{ and }\qquad}

\newcommand{\defref}[1]{Definition~\ref{#1}}

\theoremheaderfont{\fontfamily{pzc}\bfseries\large}

\newtheorem{innercustomdef}{Definition}
\newenvironment{customdef}[1]
  {\renewcommand\theinnercustomdef{#1}\innercustomdef}
  {\endinnercustomdef}
  
\newtheorem{innercustomlem}{Lemma}
\newenvironment{customlem}[1]
  {\renewcommand\theinnercustomlem{#1}\innercustomlem}
  {\endinnercustomlem}

{\theorembodyfont{\rmfamily}
\newenvironment{proof}{{\flushleft \emph{Proof}:}}{\hfill\ding{110}}

\newcommand{\g}{\mathfrak{g}}
\newcommand{\h}{\mathfrak{h}}
\newcommand{\euc}{\mathfrak{e}}
\newcommand{\Vol}{d\text{Vol}}
\newcommand{\Volume}{d\text{Vol}_\g}
\newcommand{\Volumen}{d\text{Vol}_{\g_n}}

\newcommand{\Volumeh}{d\text{Vol}_{\mathfrak{h}}}
\newcommand{\textVol}{\text{Vol}}

\newcommand{\M}{{\calM}}
\newcommand{\N}{\mathcal{N}}

\newcommand{\id}{{\text{Id}}}

\newcommand{\dist}{\operatorname{dist}}

\newcommand{\SO}[1]{\text{SO}(#1)}

\newcommand{\limn}{\lim_{n\to\infty}}

\newcommand{\tH}{\tilde{H}}

\newcommand{\e}{\varepsilon}

\newcommand{\dis}{\operatorname{dis}}
\newcommand{\Lip}{\operatorname{Lip}}

\newcommand{\Cy}[1]{{\color{red} #1}}

\renewcommand{\Cy}[1]{{ #1}}

\numberwithin{equation}{section}

\begin{document}

\title{The emergence of torsion in the continuum limit of distributed edge-dislocations -- erratum}
\author{Raz Kupferman and Cy Maor}
\date{}
\maketitle

In \cite{KM15}, following an example of locally flat Riemannian manifolds with edge-dislocation like singularities that converge to a Weitzenb\"ock manifold (Section 3), we defined a general notion of convergence of Weitzenb\"ock manifolds (Definition 4.1).
This definition had to be weak enough such that it applies to the example in Section 3, and strong enough to be well defined, that is, strong enough to allow us to prove that the limit is unique. 
The uniqueness result is Theorem 4.2 in the paper, and its proof is the main part of Section 4.

Definition 4.1 is the following:
\begin{customdef}{4.1}
\label{df:convergence}
Let $(\calM_n,\g_n,\nabla_n)$, $(\calM,\g,\nabla)$ be compact \Cy{oriented} $d$-dimensional Weitzenb\"ock manifolds with corners.
We say that the sequence $(\calM_n,\g_n,\nabla_n)$ converges to $(\calM,\g,\nabla)$ 
with $p\in \Cy{[}d, \infty)$, if there exists a sequence of diffeomorphisms   $F_n: A_n\subset \M\to \M_n$ such that:
\begin{enumerate}
\item $A_n$ covers $\M$ asymptotically:
\[
\limn \textVol_\g (\M\setminus A_n ) = 0.
\]
\item $F_n$ are approximate isometries: the distortion vanishes asymptotically, namely,
\[
\limn\dis F_n = 0.
\]
\item $F_n$ are asymptotically rigid in the mean:
\[
\limn \int_{A_n} \dist{^p}(dF_n,\SO{\g,\g_n}) \,\Volume = 0.
\]
\item The parallel transport converges in the mean in the following sense: every point in $\M$ has a neighborhood $U\subset\M$, with (i) a $\nabla$-parallel frame field $E$ on $U$, and (ii) a sequence of $\nabla_n$-parallel frame fields $E_n$  on $F_n(U\cap A_n)$,
such that
\[
\limn \int_{U\cap A_n} |F_n^\star E_n-E|^p_\g \Volume = 0.
\]
\end{enumerate}
\end{customdef}

It turns out that there is an error in the proof of Lemma 4.7 (Lemma~4.8 in the arXiv version of \cite{KM15}), which is a part of the proof of the uniqueness of limit (Theorem 4.2). 
In order to overcome it, one has to strengthen the assumptions in \defref{df:convergence}.
A simple way of doing so is by demanding that there exists a constant $C$ such that $\Lip(F_n), \Lip(F_n^{-1}) < C$ for every $n$, that is, by assuming that $F_n$ are uniformly bi-Lipschitz. 
This makes the proof significantly simpler (in particular, the widely used Lemma~4.6 becomes trivial as the sets $A_n^\e$ are eventually equal to $A_n$). 
This assumption also makes the requirement $p\ge d$ irrelevant -- if the assumptions hold for any $p\ge 1$, then they hold for any $p<\infty$. This is explained in detail in \cite{KM16}, a related paper that makes this bi-Lipschitz assumption (for other reasons).

While the uniform bi-Lipschitz assumption is a restrictive assumption, the example presented in \cite{KM15}, as well as the general construction of the same phenomenon presented in \cite{KM15b} all involve uniformly bi-Lipschitz mappings.

However, we find this assumption a bit unnatural and too restrictive, and so we prefer to present here an intermediate one, stronger than \defref{df:convergence} but weaker than \defref{df:convergence} + the uniform bi-Lipschitz assumption:
\begin{customdef}{4.1'}
\label{df:convergence_new}
Let $(\calM_n,\g_n,\nabla_n)$, $(\calM,\g,\nabla)$ be compact \Cy{oriented} $d$-dimensional Weitzenb\"ock manifolds with corners.
Let $p_{\text{min}}(d) = d + 1 + \frac{d}{2}(\sqrt{1+4/d^2}-1)$.
We say that the sequence $(\calM_n,\g_n,\nabla_n)$ converges to $(\calM,\g,\nabla)$ 
with $p\in [p_{\text{min}}, \infty)$, if there exists a sequence of diffeomorphisms   $F_n: A_n\subset \M\to \M_n$ such that:
\begin{enumerate}
\item $A_n$ covers $\M$ asymptotically:
\[
\limn \textVol_\g (\M\setminus A_n )\, \Lip(F_n)^2 = \limn \textVol_\g (\M\setminus A_n )\, \Lip(F_n^{-1})^2 = 0.
\]
\item $F_n$ are approximate isometries: the distortion vanishes asymptotically, namely,
\[
\limn\dis F_n = 0.
\]
\item $F_n$ are asymptotically rigid in the mean:
\[
\begin{split}
&\limn \int_{A_n} \dist{^p}(dF_n,\SO{\g,\g_n}) \,\Volume =0, \\
&\limn \int_{\M_n} \dist{^p}(dF_n^{-1},\SO{\g_n,\g}) \,\Volumen = 0.
\end{split}
\]
\item The parallel transport converges in the mean in the following sense: every point in $\M$ has a neighborhood $U\subset\M$, with (i) a $\nabla$-parallel frame field $E$ on $U$, and (ii) a sequence of $\nabla_n$-parallel frame fields $E_n$  on $F_n(U\cap A_n)$,
such that
\[
\limn \int_{U\cap A_n} |F_n^\star E_n-E|^p_\g \Volume = 0.
\]
\end{enumerate}
\end{customdef}

The differences between \defref{df:convergence} and \defref{df:convergence_new} are:
\begin{enumerate}
\item The condition on $p$ is more restrictive (instead of $p\ge d$ we assume $p\ge p_{\text{min}}(d)$, where  $p_{\text{min}}^2 - (d+2) p_{\text{min}} +d =0$).
\item Condition (1) now relates the size of the "holes" in $\M$ to the "wildness" of $F_n$.
\item Condition (3) now requires that $F_n$ and $F_n^{-1}$ are both asymptotically rigid. 
	That is, there is a symmetric penalization for both expansion and contraction, instead of penalizing mainly expansions. 
	Adding a penalization for large contractions is very natural from the material science and elasticity point of view, which is the main motivation for this work.
\end{enumerate}

Below is a restatement of Lemma 4.7 in \cite{KM15} (Lemma~4.8 in the arXiv version), and a proof under the assumptions of \defref{df:convergence_new}, which shows that the limit is indeed unique as stated in Theorem~4.2.
\begin{customlem}{4.7'}
\label{lm:uniqueness}
Let $(\M_n,\g_n)$, $(\M,\g)$ and $(\N,\mathfrak{h})$ be compact Riemannian manifolds.
Let $E_n$, $E^\M$ and $E^\N$ be frame fields on $\M_n$, $\M$ and $\N$, respectively.
Suppose that both 
\[
(\M_n,\g_n,E_n) \to (\M,\g,E^\M)
\Textand
(\M_n,\g_n,E_n) \to (\N,\h,E^\N)
\]
with respect to diffeomorphisms $F_n:A_n\subset\M \to \M_n$ and $G_n:B_n\subset\N\to \M_n$ (here, the pullbacks of the frame fields converge in $L^p$).
Then $H_\star E^\M = E^\N$, where $H:\M\to\N$ is the uniform limit of $H_n = G_n^{-1}\circ F_n$  
defined in Lemma~4.3. 
\end{customlem}

\begin{proof}
We need to show that $H_\star E^\M - E^N = 0$. Since $H$ is the limit of $H_n$, we start by 
estimating $(H_n)_\star E^\M - E^\N$.
We fix some $\e>0$, and consider $H_n$ as a diffeomorphism $A_n^\e\to H_n(A_n^\e)$, where sets $A_n^\e$ are defined in Lemma~4.6.
By the standard inequality $|a+b|^p \le C(|a|^p + |b|^p)$ we get
\[
\begin{split}
&\int_{A_n^\e} |dH_n E^\M - H_n^*E^\N|_{H_n^*\h}^p \Volume \\
	&\quad \le C\int_{A_n^\e} |dH_n E^\M - H_n^*G_n^\star E_n|_{H_n^*\h}^p \Volume
		+ C \int_{A_n^\e} | H_n^* G_n^\star E_n - H_n^* E^\N|_{H_n^*\h}^p \Volume \\
	&\quad = C\int_{A_n^\e} |dH_n E^\M - H_n^*G_n^\star E_n|_{H_n^*\h}^p \Volume
		+ C \int_{H_n(A_n^\e)} | G_n^\star E_n - E^\N|_{\h}^p \frac{\Vol_{(H_n)_\star\g}}{\Volumeh} \Volumeh \\
	&\quad \le C' \int_{A_n^\e} |dH_n|^p |E^\M - F_n^\star E_n|_{H_n^\star\h}^p \Volume
		+ C \int_{H_n(A_n^\e)} | G_n^\star E_n - E^\N|_{\h}^p \frac{\Vol_{(H_n)_\star\g}}{\Volumeh} \Volumeh \\
	&\quad \le C'' \int_{A_n^\e} |E^\M - F_n^\star E_n|_{\g}^p \Volume
		+ C'' \int_{H_n(A_n^\e)} | G_n^\star E_n - E^\N|_{\h}^p \Volumeh,
\end{split}
\]
where we used the uniform bounds on $|dH_n|$ and $|dH_n^{-1}|$ on $A_n^\e$, and Lemma~4.5.
Now, the first addend in the last line tends to $0$ since $(\M_n,\g_n,E_n) \to (\M,\g,E^\M)$ with respect to the maps $F_n$, and the second addend since $(\M_n,\g_n,E_n) \to (\N,\h,E^\N)$ with respect to the maps $G_n$. 
Therefore, we have established that
\beq
\label{eq:tmp1}
\int_{A_n^\e} |dH_n E^\M - H_n^*E^\N|_{H_n^*\h}^p \Volume  \to 0.
\eeq
The proof would be complete if we could replace $(H_n)_\star$ by $H_\star$ and $H_n(A_n^\e)$ by $\N$ in the limit $n\to\infty$. This is not yet possible since $H_n$ tends to $H$ on $A_n$ only uniformly, whereas the push-forward of frame fields with $H_n$ involves derivatives of $H_n$. 
Therefore, we will show that $dH_n\to dH$ in some sense. 

We start by showing that
\beq
\label{eq:H_n_to_SO}
\limn \int_{A_n} \dist(dH_n,\SO{\g,\h}) \,\Volume = 0.
\eeq
Indeed, let $x\in A_n$, and let $q_n\in \SO{\g_x,(\g_n)_{F_n(x)}}$ be a point that realizes $\dist(dF_n,\SO{\g,\g_n})$ at $x$, and $r_n\in  \SO{(\g_n)_{F_n(x)},\h_{H_n(x)}}$ a point that realizes $\dist(dG_n^{-1},\SO{\g_n,\h})$ at $F_n(x)$. 
Then we have at the point $x$,
\[
\begin{split}
\dist(dH_n,\SO{\g,\h}) &\le |dH_n - r_n q_n| = |dG_n^{-1}dF_n - r_n q_n| \\
	& \le |dG_n^{-1} dF_n - r_n dF_n| + |r_n dF_n -r_n q_n| \\
	& \le |dG_n^{-1} - r_n | |dF_n| + |dF_n - q_n| \\
	& = \dist(dG_n^{-1},\SO{\g_n,\h})|dF_n| + \dist(dF_n,\SO{\g,\g_n})
\end{split}
\]
and therefore globally
\[
\dist(dH_n,\SO{\g,\h}) \le F_n^*\dist(dG_n^{-1},\SO{\g_n,\h})\,|dF_n| + \dist(dF_n,\SO{\g,\g_n})
\]
The second addend vanishes in $L^p(A_n,\g)$ as $n\to \infty$ and therefore also in $L^1(A_n,\g)$.
As for the first addend, using H\"older inequality and Lemma~4.5, we obtain:
\[
\begin{split}
& \| F_n^*\dist(dG_n^{-1},\SO{\g_n,\h})\,|dF_n|\|_{L^1(A_n,\g)} \\
	&\quad \le \| F_n^*\dist^{p/p-1}(dG_n^{-1},\SO{\g_n,\h})\|_{L^{1}((\M,\g))}^{(p-1)/p} \,\| dF_n\|_{L^p(A_n,\g)} \\
	&\quad = \| \dist^{p/p-1}(dG_n^{-1},\SO{\g_n,\h}) \,\frac{\Vol_{(F_n)_\star\g}}{\Volumen} \|_{L^{1}(\M_n,\g_n)}^{(p-1)/p} \,\| dF_n\|_{L^p((A_n,\g))} \\
	&\quad \le \| \dist^{p/p-1}(dG_n^{-1},\SO{\g_n,\h}) \,|dF_n^{-1}|^d \|_{L^{1}(\M_n,\g_n)}^{(p-1)/p} \,\| dF_n\|_{L^p(A_n,\g)} \\
	&\quad \le \| \dist^{p/p-1}(dG_n^{-1},\SO{\g_n,\h})\|_{L^{p/p-d}(\M_n,\g_n)}^{(p-1)/p}\, \| |dF_n^{-1}|^d\|_{L^{p/d}((\M_n,\g_n))}^{(p-1)/p}\, \|dF_n\|_{L^p(A_n,\g)} \\
	&\quad = \| \dist(dG_n^{-1},\SO{\g_n,\h})\|_{L^{p^2/(p-d)(p-1)}(\M_n,\g_n)}^{\alpha}\, \| dF_n^{-1}\|_{L^{p}(\M_n,\g_n)}^{\beta}\, \|dF_n\|_{L^p(A_n,\g)},
\end{split}
\]
where $\alpha,\beta>0$ are the appropriate powers (they are immaterial for the rest of the argument).
Now, the last two terms on the last line are uniformly bounded in $n$ by our assumptions on $dF_n$.
The first term vanishes as $n$ goes to infinity by our assumptions on $G_n$, since our assumption on $p$ implies (i) $p^2/(p-d)(p-1)\le p$, and (ii) $p>d$, hence $\textVol_{\g_n}(\M_n)\to \textVol_{\g}(\M)$ (this is an immediate corollary of Lemma~4.5) and so the constants in H\"older inequality used to replace $L^{p^2/(p-d)(p-1)}(\M_n,\g_n)$ with $L^{p}(\M_n,\g_n)$ are bounded uniformly in $n$. 

Recall that we want to prove convergence $dH_n\to dH$ in an appropriate sense. 
Since Sobolev spaces are easier to handle when the image is a vector bundle,
we fix an isometric immersion $\phi:(\N,\h)\to (\R^\nu,\euc)$ for  large enough $\nu$, where $\euc$ is the standard Euclidean metric. 
Consider the mappings $\phi\circ H_n:A_n\to \R^\nu$. 
These mappings satisfy
	\beq
	\label{eq:phiH_n_to_O}
	\limn \int_{A_n} \dist(d(\phi\circ H_n),\text{O}(\g,\euc)) \,\Volume = 0,
	\eeq
since the left hand side is bounded from above by the left hand side of \eqref{eq:H_n_to_SO}, and $\phi$ is an isometric immersion.

Since $\phi\circ H_n$ are smooth, and in particular Lipschitz, we can extend them to $\bar{H}_n\in W^{1,\infty}(\M;\R^\nu)$, such that $\|d\bar{H}_n\|_{\infty}= \Lip(\bar{H}_n) \le C\Lip(\phi\circ H_n)$ for some constant independent of $n$ (for example, we can use McShane extension lemma \cite{Hei05}, or more sophisticated results with a better constant $C$).
We claim that $\bar{H}_n$ satisfy
\beq
\label{eq:barH_n_to_O}
\limn \int_{\M} \dist(d\bar{H}_n,\text{O}(\g,\euc)) \,\Volume = 0.
\eeq
Indeed,
\[
\begin{split}
& \int_{\M} \dist(d\bar{H}_n,\text{O}(\g,\euc)) \,\Volume \\
	&\quad \le \int_{A_n} \dist(d\bar{H}_n,\text{O}(\g,\euc)) \,\Volume + (C_1+\|d\bar{H}_n\|_\infty)\textVol_\g(\M\setminus A_n) \\
	&\quad \le \int_{A_n} \dist(d(\phi \circ H_n),\text{O}(\g,\euc)) \,\Volume + C_1\textVol_\g(\M\setminus A_n) \\
		&\quad\quad + C_2\Lip(\phi\circ H_n)\textVol_\g(\M\setminus A_n) 
\\
	&\quad \le \int_{A_n} \dist(d(\phi \circ H_n),\text{O}(\g,\euc)) \,\Volume +  C_1\textVol_\g(\M\setminus A_n) \\
		&\quad\quad + C_3\Lip(G_n^{-1}) \Lip(F_n) \textVol_\g(\M\setminus A_n) \\
	&\quad = \int_{A_n} \dist(d(\phi \circ H_n),\text{O}(\g,\euc)) \,\Volume  +  C_1\textVol_\g(\M\setminus A_n)  \\
	&\quad \quad + C_3\sqrt{\Lip(G_n^{-1})^2 \textVol_\g(\M\setminus A_n)} \sqrt{ \Lip(F_n)^2 \textVol_\g(\M\setminus A_n)}.
\end{split}
\]
The first summand goes to zero by \eqref{eq:phiH_n_to_O}, and the first summand by the asymptotic surjectivity of $F_n$, which also imply that $ \Lip(F_n)^2 \textVol_\g(\M\setminus A_n)\to 0$.
We are left to deal with the term $\Lip(G_n^{-1})^2 \textVol_\g(\M\setminus A_n)$.
Note that by moving to a subsequence, we can assume that $ \textVol_\g(\M\setminus A_n)/ \textVol_\h(\N\setminus B_n)$ is monotone. 
Since the roles of $\M$ and $\N$ (and their associated metrics, mappings, etc.) are completely symmetric, we can assume without loss of generality that this sequence is monotonically decreasing, and in particular, bounded, hence 
\[
\Lip(G_n^{-1})^2 \textVol_\g(\M\setminus A_n) \le C\Lip(G_n^{-1})^2 \textVol_\g(\N\setminus B_n) \to 0
\]
by our assumptions on $G_n$.

We therefore establish \eqref{eq:barH_n_to_O}. 
While $\bar{H}_n$ are Lipschitz functions on $\M$, they are not uniformly Lipschitz.
In order to complete the proof, we will replace them by uniformly Lipschitz maps $\tH_n:\M\to \R^\nu$ that agree with $\bar{H}_n$ over large sets.
That is, we now claim that there exist maps $\tH_n:\M\to \R^\nu$ such that
\begin{enumerate}
\item $\textVol_{\g}\{\tH_n \ne \bar{H}_n\} \to 0 $, and
\item $\tH_n$ are uniformly bounded in $W^{1,\infty}(\M;\R^\nu)$, and in particular Lipschitz by with a uniform constant $L$.
\end{enumerate}
To show this, we use Proposition A.1 in \cite{FJM02b} on $\bar{H}_n$, with $\lambda$ large enough such that for every $x$ and every $T:T_x\M\to \R^\nu$
\[
|T|\ge \lambda \quad \Rightarrow |T| < 2\dist(T,\text{O}(\g,\euc)).
\]
Proposition A.1 in \cite{FJM02b} then implies that there exists a sequence of functions $\tH_n$, uniformly bounded in $W^{1,\infty}(\M;\R^\nu)$ such that
\[
\textVol_{\g}(\M\setminus R_n)  \le C\int_{\{x\in\M : |d \bar{H}_n|>\lambda \}} |d \bar{H}_n| \,\Vol_{\g},
\] 
where $R_n := \{\tH_n = \bar{H}_n\}$. 
Therefore, we have
\[
\begin{split}
\textVol_{\g}(\M\setminus R_n) & \le C\int_{\{x\in\M : |d \bar{H}_n|>\lambda \}} |d \bar{H}_n| \,\Vol_{\g} \\
	& \le 2C\int_{\{x\in\M : |d \bar{H}_n|>\lambda \}}\dist(d\bar{H}_n,\text{O}(\g,\euc)) \,\Vol_{\g} \\
	& \le 2C\int_{\M}\dist(d\bar{H}_n,\text{O}(\g,\euc)) \,\Vol_{\g} \to 0.
\end{split}
\]
This argument is similar to Lemma~3.3 in \cite{LP10}.\footnote{Note that while Proposition A.1 in \cite{FJM02b} discusses a Lipschitz domain in the Euclidean space, and therefore directly applies to a manifold that can be covered by a single coordinate chart (as in \cite{LP10}), this is not a problem here:
First, the claim of Lemma \ref{lm:uniqueness} is local, hence we could work locally and assume w.l.o.g.~that $\M$ is covered by a single chart.
Second, looking more carefully at the proof of Proposition A.1 in \cite{FJM02b}, the same partition of unity argument used there to discuss the general Lipschitz domain can actually be used again to discuss a general Riemannian manifold.}

The functions $\tH_n$ converge to $\phi\circ H$ uniformly on $\M$, as 
\[
\begin{split}
d_{\R^\nu}(\tH_n(p), \phi\circ H(p)) &\le 
	d_{\R^\nu}(\tH_n(p),\tH_n(\psi_n(p))) + d_{\R^\nu}(\tH_n(\psi_n(p)),\phi\circ H(\psi_n(p))) \\
&\qquad+ d_{\R^\nu}(\phi\circ H(\psi_n(p)),\phi\circ H(p))\\
	&= d_{\R^\nu}(\tH_n(p),\tH_n(\psi_n(p))) + d_{\R^\nu}(\phi\circ H_n(\psi_n(p)),\phi\circ H(\psi_n(p))) \\
&\qquad+ d_{\R^\nu}(\phi\circ H(\psi_n(p)),\phi\circ H(p))\\
&\le d_{\R^\nu}(\tH_n(p),\tH_n(\psi_n(p))) + d_\N(H_n(\psi_n(p)),H(\psi_n(p))) \\
&\qquad+ d_\N(H(\psi_n(p)),H(p))\\
&=  d_{\R^\nu}(\tH_n(p),\tH_n(\psi_n(p))) + d_\N(H_n(\psi_n(p)),H(\psi_n(p))) \\
&\qquad+ d_\M(\psi_n(p),p) \\
	&\le L\cdot d_\M(p,\psi_n(p)) + d_\N(H_n(\psi_n(p)),H(\psi_n(p))) + d_\M(\psi_n(p),p) \\
	&\le (L+1)\,\sup_{\M} d(\cdot,\psi_n(\cdot)) + \sup_{A_n} d_\N(H_n(\cdot),H(\cdot)) \to 0.
\end{split}
\]
Here $\psi_n$ is a mapping $\M\to A_n\cap R_n$ satisfying
\[
\psi_n|_{A_n\cap R_n} = \id
\Textand
\sup_{p\in\M} d_\M(p,\psi_n(p)) < \e_n
\]
for some $\e_n\to 0$; it is analogous to the mapping $\M\to A_n$ introduced in Lemma~4.3.
Since $\textVol_\g(\M\setminus R_n)\to 0$, we can indeed choose such a sequence $\e_n\to 0$.
In the passage from the first to the second line we used the fact that $\tH_n$ coincides with $\phi\circ H_n$ on the image of $\psi_n$.  In the passage from the second to the third line we used the fact that $\phi$ is distance reducing. In the passage from the third to the fourth line we used the fact that $H$ is an isometry. The rest follows from the uniform Lipschitz bound on $\tH_n$ and the uniform convergence of $\psi_n$ to $\id_\M$, and  the uniform convergence of $H_n$ to $H$ on $A_n$.

Now, note that \eqref{eq:tmp1} can be written as
\beq
\int_{A_n^\e} |d(\phi\circ H_n) E^\M - H_n^*(d\phi\circ E^\N)|_{\euc}^p \Volume \to 0.
\eeq
Therefore,
\[
\begin{split}
& \int_\M |d\tH_n\circ E^\M - \tH_n^*\phi_\star E^\N|_{\euc}^p \Vol_{\g} \\
	& \quad \le \int_{A_n^\e\cap R_n} |d\tH_n\circ E^\M - \tH_n^*\phi_\star E^\N|_{\euc}^p \Vol_{\g} \\
		& \quad \quad + C\,\|d\tH_n\circ E^\M - \tH_n^*\phi_\star E^\N\|_\infty^p\textVol_\g(\M\setminus(A_n^\e\cap R_n)) 
\end{split}
\]
\[
\begin{split}
	& \quad = \int_{A_n^\e\cap R_n} |d(\phi\circ H_n) E^\M - H_n^*(d\phi\circ E^\N)|_{\euc}^p \Vol_{\g} \\
		& \quad \quad + C\,\|d\tH_n\circ E^\M - \tH_n^*\phi_\star E^\N\|_\infty^p\textVol_\g(\M\setminus(A_n^\e\cap R_n)) \\
	& \quad \le \int_{A_n^\e} |d(\phi\circ H_n) E^\M - H_n^*(d\phi\circ E^\N)|_{\euc}^p \Vol_{\g} \\
		&\quad \quad + C\,(L\| E^\M\|_\infty + \|E^\N\|_\infty)^p\textVol_\g(\M\setminus(A_n^\e\cap R_n)) \to 0,
\end{split}
\]
where we used the the uniform Lipschitz constant of $H_n$ (denoted by $L$), the fact that $\textVol_\g(\M\setminus A_n^\e)\to 0$ by Lemma~4.6 and $\textVol_\g(\M\setminus R_n)\to 0$.

Since $\tH_n\to \phi\circ H$ uniformly and $E^\N$ is smooth, we can replace $\tH_n^*$ by $(\phi\circ H)^*$. Therefore we obtain
\[
\int_{\M} |d\tH_n\circ E^\M - H^*(d\phi\circ E^\N)|_{\euc}^p \Vol_{\g} \to 0.
\]
It follows that $d\tH_n$ converges in $L^p(\M;T^*\M\otimes\R^\nu)$ to the map
\[
E^\M \mapsto H^*(d\phi\circ E^\N).
\]
Since, in addition, $\tH_n$ converges uniformly to $\phi\circ H$, it follows that $\tH_n$ converges to $\phi\circ H$ in $W^{1,p}(\M;\R^\nu)$, and in particular,
\[
d(\phi\circ H) \circ E^\M = H^*(d\phi\circ E^\N).
\]
Since $\phi$ is an embedding we can eliminate $d\phi$ on both sides, getting
\[
H_\star E^\M = E^\N.
\]

\end{proof}

{\footnotesize
\bibliographystyle{amsalpha}

\providecommand{\href}[2]{#2}
\providecommand{\arxiv}[1]{\href{http://arxiv.org/abs/#1}{arXiv:#1}}
\providecommand{\url}[1]{\texttt{#1}}
\providecommand{\urlprefix}{URL }
}

\end{document}